\documentclass{article}
\usepackage{amsmath,amssymb,amsthm}
\usepackage{amsfonts}
\usepackage{mathrsfs}

\newtheorem{theorem}{Theorem}[section]
\newtheorem{definition}{Definition}[section]
\newtheorem{corollary}{Corollary}[section]
\newtheorem{question}{Question}[section]
\newtheorem{remark}{Remark}[section]
\newtheorem{proposition}{Proposition}[section]

\usepackage{geometry}
\title{On strongly convex K\"ahler--Berwald manifolds}
\author{
Xiaoshu Ge \thanks{An early version of this manuscript was completed in December 2022, when the author was a PhD student at the School of Mathematical Sciences, Xiamen University.}
\\
\small Institute of Mathematics, Academy of Mathematics and Systems Science, \\ \small Chinese Academy of Sciences, Beijing, 100190, China\\ \small
E-mail:gxs@amss.ac.cn
}
\date{}

\begin{document}
\maketitle

\begin{abstract}
%Deng proved that every real symmetric Finsler metric is a Berwald metric, for an irreducible Hermitian symmetric space with $G$-invariant Finsler metric $F$ there exists an invariant complex structure $J$ such that $(G / H, J, F)$ is a symmetric complex Finsler metric on $G / H$. %Thus one may wonder whether $(G / H, J, F)$ is a K\"ahler--Berwald metric?

For every real Berwald strongly convex complex Finsler manifold $(M,F,J,\hat\nabla)$, if the Cartan nonlinear connection satisfies the parallelism condition $\hat\nabla J=0$,
we construct a good complex vertical connection $D$,
prove that $D$ is the complex Rund connection, and $F$ is K\"ahler-Berwald. The inverse direction also can be seen from the construction and proof.
As an application, we prove that on an irreducible Hermitian symmetric space $M=G/H$, $G$-invariant complex Finsler metrics are K\"ahler--Berwald.
Besides, we obtain the characterization of strongly pseudoconvex holomorphically invariant complex Finsler metrics on the classical domains of type-IV.
\end{abstract}


\section{Introduction and main results}
The K\"ahler property can be characterized by the parallelism of the almost complex structure with respect to the Riemannian connection. It is still a question for Finsler geometry before the manuscript written in December 2022:
\begin{question}\label{ques1}
For a strongly convex complex Finsler metric, wether the complex structure is parallel with respect to the nonlinear connection associated to the Cartan connection coincides with $F$ being K\"ahler--Finsler?
\end{question}

In \cite{Zh11}, Zhong proved that for a strongly convex weakly K\"ahler Finsler metric, it is a weakly complex Berwald metric if and only if it is a real Berwald metric; and in \cite{Xia-Z}, Xia and Zhong proved that for a complex Berwald metric $F$, it is a K\"ahler Finsler metric if and only if it is a weakly K\"ahler Finsler metric. Thus, to obtain a strongly convex K\"ahler Berwald metric, it is nature to assume $F$ is a real Berwald metric.
In this paper, we consider real Berwald strongly convex complex Finsler manifold $(M,F,J,\hat\nabla)$.
If the Cartan nonlinear connection satisfies the parallelism condition $\hat\nabla J=0$, we construct a good complex vertical connection $D$, prove that $D$ is the complex Rund connection, and $F$ is K\"ahler-Berwald.
\begin{theorem} \label{thm}
Let $(M, J, F, \hat{\nabla})$ be a real Berwald strongly convex complex Finsler manifold. If $\hat \nabla J=0$, then $D: \Gamma(\mathcal V) \to \Gamma(T^{*} \tilde{M} \otimes \mathcal V)$,
\[D(Y):=\left(d Y^{\alpha}+Y^{\beta} \Gamma_{\beta ; \gamma}^{\alpha} d z^{\gamma}\right) \otimes \dot{\partial}_{\alpha},\quad \forall Y=Y^{\alpha} \dot{\partial}_{\alpha} \in \Gamma(\mathcal{V}), \] is the complex Rund connection, and $F$ is K\"ahler--Berwald.
\end{theorem}
\begin{remark}
The inverse direction also can be seen from the construction and proof.
\end{remark}

In \cite{Deng}, S. Deng proved that every real symmetric Finsler metric is a Berwald metric, for an irreducible Hermitian symmetric space with $G$-invariant Finsler metric $F$ there exists an invariant complex structure $J$ such that $(G / H, J, F)$ is a symmetric complex Finsler metric on $G / H$. Thus one may wonder whether $(G / H, J, F)$ is a K\"ahler--Berwald metric?
As a corollary of Theorem \ref{thm}, we obtain
\begin{theorem}\label{Symm}
On an irreducible Hermitian symmetric space $(M=G / H, J)$, any $G$-invariant strongly convex complex Finsler metric $F$ is a K\"ahler Berwald metric.
\end{theorem}

At last, we prove that all holomorphically invariant strongly pseudoconvex complex Finsler metrics of the classical domains of type IV are given by a class of examples provided in \cite{GZ}.
\begin{theorem}\label{R4}
For every holomorphically invariant strongly pseudoconvex complex Finsler metric $F$ of $\mathcal R_{IV}$, its isotropic invariant norm at the origin can be written as follows
\begin{eqnarray}
	F(0;v)=:f(v)=\sqrt{r\phi(s)},\quad r=v\bar{v}', s(v)=\frac{|vv'|^2}{r^2}, 0\neq v\in T_{0}^{1,0} \mathfrak{R}_{IV}.
\end{eqnarray}
\end{theorem}

\begin{remark}
Theorems \ref{thm}--\ref{Symm} are proved in Sections \ref{sec2}--\ref{sec3}, and the proof of Theorem \ref{R4} is given in Section \ref{R4}.

An early version of this manuscript was completed in December 2022, with part of the material for Theorem \ref{R4} finished as early as October 2022, and we make no change about any proof of that manuscript. The earlier unpublished manuscript was shared with and reviewed by Professor Chunping Zhong in 2022.

Several related works appeared subsequently, all after the completion of our 2022 draft:
\begin{itemize}
  \item In \cite{Zho}, C. Zhong adopts the same core idea for the proof of Theorem \ref{R4} as that developed in this manuscript.
  \item In \cite{XZ}, W. Xiao and C. Zhong adopt the identical construction of the complex Rund connection from the Cartan connection, as well as the proof methods used in Theorem \ref{thm} of the present work. They construct the Chern--Finsler connection, whose horizontal part coincides exactly with the complex Rund connection introduced here, and obtain the equivalence that $\hat\nabla J=0$ if and only if $F$ is K\"ahler--Berwald, in parallel with our result.
  \item In \cite{HZ}, Y. He and C. Zhong establish results essentially identical to Corollary \ref{Symm}. Their argument also relies on the fact that $\nabla J=0$, combined with a real Berwald-type result from \cite{Deng}, to conclude that $F$ is K\"ahler--Berwald, following the same line of reasoning as our work.
\end{itemize}

We explicitly claim full priority for all original geometric constructions, proof ideas and main results presented in this paper.
\end{remark}

\section{Strongly convex complex Finsler metrics}\label{sec2}
A strongly convex complex Finsler metric is a complex Finsler metric which also canonically defines a real Finsler metric as we shall introduce below. Thus let us recall some basic notions of real and complex Finsler geometry at first.

Let $(M, J)$ be a complex manifold of dimension $n\,(n \geq 2)$ with complex structure $J$, $TM$ be the real tangent bundle of $M$, and $z=(z_{1}, \cdots, z_{n})$ be local complex coordinates with $z_{j}=x_{j}+\sqrt{-1}\, x_{j+n}$. Then
\[
J\left( \frac{\partial }{\partial x_{j}}\right) =\frac{\partial }{\partial x_{j+n}},\quad J\left( \frac{\partial }{\partial x_{j+n}}\right) =-\frac{\partial }{\partial x_{j}},\quad j=1,\cdots ,n.
\]

Extend $J$ by complex linearity into the complexified tangent bundle $T^{\mathbb{C}} M$ of $M$. Then we have the direct sum $T_{z}^{\mathbb{C}} M=T_{z}^{1,0} M \oplus T_{z}^{0,1} M$, where $T_{z}^{1,0} M$ and $T_{z}^{0,1} M$ are called (complexified) tangent vectors of type $(1, 0)$ and $(0, 1)$, respectively, as the $\pm\sqrt{-1}$-eigenspaces of $J$. We also have the bundle decomposition $T^{\mathbb{C}} M=T^{1,0} M \oplus T^{0,1} M$. A basis of $T_{z}^{1,0} M$ is given by $\frac{\partial}{\partial z_{1}}, \cdots, \frac{\partial}{\partial z_{n}}$, where
\[
\frac{\partial}{\partial z_{j}}=\frac{1}{2}\left(\frac{\partial}{\partial x_{j}}-\sqrt{-1}\, \frac{\partial}{\partial x_{j+n}}\right).
\]
Similarly, a basis of $T_{z}^{0,1} M$ is given by $\frac{\partial}{\partial \bar{z}_{1}}, \cdots, \frac{\partial}{\partial \bar{z}_{n}}$, where
\[
\frac{\partial}{\partial \bar{z}_{j}}=\frac{1}{2}\left(\frac{\partial}{\partial x_{j}}+\sqrt{-1}\, \frac{\partial}{\partial x_{j+n}}\right).
\]

In the following, we shall denote by $\tilde{M}$ the complement of the zero section in $T^{1,0} M$, and use a semi-colon to distinguish derivatives with respect to base manifold coordinates $z=(z^{1}, \cdots, z^{d})$ and fiber coordinates $v=(v^{1}, \cdots, v^{d})$. We use Greek indices and Latin indices to distinguish derivatives with respect to complex coordinates $z,v$ and real coordinates $x,u$, respectively. For example,
\[
G_{\alpha }=\frac{\partial G}{\partial v^{\alpha }},\quad G_{\alpha \overline{\beta}}=\frac{\partial ^{2}G}{\partial v^{\alpha }\partial \overline{v}^{\beta }},\quad G_{;\alpha }=\frac{\partial G}{\partial z^{\alpha }},\quad G_{\alpha ;\overline{\beta}}=\frac{\partial ^{2}G}{\partial v^{\alpha }\partial \overline{z}^{\beta }};
\]
\[
G_{i}=\frac{\partial G}{\partial u^{i}},\quad G_{ij}=\frac{\partial ^{2}G}{\partial u^{i}\partial u^{j}},\quad G_{;i}=\frac{\partial G}{\partial x^{i}},\quad G_{i;j}=\frac{\partial ^{2}G}{\partial u^{i}\partial x^{j}}.
\]

For simplicity we set
\[
\partial _{i}:=\frac{\partial }{\partial z^{i}},\quad \partial _{\overline{i}}:=\frac{\partial }{\partial \overline{z}^{i}},\quad \dot{\partial}_{i}:=\frac{\partial }{\partial v^{i}},\quad \dot{\partial}_{\overline{i}}:=\frac{\partial }{\partial \overline{v}^{i}};\quad \partial _{i}^{o}:=\frac{\partial }{\partial x^{i}},\quad \dot{\partial}_{i}^{o}:=\frac{\partial }{\partial u^{i}}.
\]

We can also define a complex structure on $T^{\mathbb{C}} \tilde{M}$, still denoted by $J$, with
\[
J\left(\partial_{i}\right)=\sqrt{-1}\, \partial_{i},\quad J\left(\partial_{\overline{i}}\right)=-\sqrt{-1}\, \partial_{\overline{i}},\quad J\left(\dot{\partial}_{i}\right)=\sqrt{-1}\, \dot{\partial}_{i},\quad J\left(\dot{\partial}_{\overline{i}}\right)=-\sqrt{-1}\, \dot{\partial}_{\overline{i}}.
\]

A complex Finsler metric on $M$ is given by \cite{AP}:
\begin{definition}
A complex (resp. real) Finsler metric on $M$ is a function $F: T^{1,0} M\,(\text{resp. }TM ) \to [0,+\infty)$ satisfying:
\begin{enumerate}
\item (Positivity) $F(v)>0$ for any non-zero $v \in T_{z}^{1,0} M\,(\text{resp. }TM )$, and equality holds iff $v=0$;
\item (Regularity) $F^{2}$ is continuous on $T^{1,0} M\,(\text{resp. }TM )$ and smooth at any $(z, v) \in \tilde{M}\,(\text{resp. }TM \setminus \{\text{zero section}\})$;
\item (Homogeneity) $F(z, \lambda v)=|\lambda| F(z, v)$ for any $(z, v) \in T^{1,0} M\,(\text{resp. }TM )$ and any $\lambda \in \mathbb{C}\,(\text{resp. } \mathbb{R})$.
\end{enumerate}
Denote $G=F^{2}$. If the $n\times n$ Hermitian matrix $(G_{\alpha \bar{\beta}})$ is positive definite for any $(z, v) \in \tilde{M}$, then $F$ is called a \emph{strongly pseudoconvex complex Finsler metric}.

Write $v=(v_{1}, \cdots, v_{n})$ with $v_{j}=u_{j}+\sqrt{-1}\, u_{j+n}$ for $j=1, \cdots, n$. If the $2n\times 2n$ matrix $(G_{ij})$ is positive definite for any $(x, u) \in TM \setminus \{\text{zero section}\}$, then $F$ is called a \emph{strongly convex complex Finsler metric}.
\end{definition}

In the following, we always assume $F$ to be a strongly convex complex Finsler metric. Since $(G_{ij})$ is positive definite for any $(x, u) \in TM \setminus \{\text{zero section}\}$, it induces a Riemannian metric $\langle \cdot | \cdot \rangle$ on the vertical subbundle $\mathcal V_{\mathbb{R}} \subset TTM$, where $\mathcal V_{\mathbb{R}}:=\ker(d \pi)$ is the kernel of the projection $\pi: TM \to M$:
\[
\langle X | Z\rangle:=\frac{1}{2} G_{ij}(u) X^{i} Z^{j},\quad \forall X, Z \in \mathcal{V}_{\mathbb{R}}.
\]

Then there exists the Cartan connection \cite{AP}:
\begin{theorem}
Let $F: TM \to \mathbb{R}^{+}$ be a Finsler metric. Then there is a unique vertical connection $\hat{D}: \Gamma(\mathcal V_{\mathbb{R}}) \to \Gamma(T^{*} \tilde{M} \otimes \mathcal V_{\mathbb{R}})$ such that:
\begin{enumerate}
\item[(i)] $\hat{D}$ is good;
\item[(ii)] for all $X \in T\tilde{M}$,
\[
X\langle V| W\rangle =\langle \hat{\nabla}_{X}V| W\rangle +\langle V| \hat{\nabla}_{X}W\rangle;
\]
\item[(iii)] $\theta(V, W)=0$ for all $V, W \in \mathcal V$, where $\theta$ is the torsion of the linear connection on $\tilde{M}$ induced by $\hat{D}$;
\item[(iv)] $\theta(H, K) \in V$ for all $H, K \in \mathcal H$.
\end{enumerate}
\end{theorem}

Similarly, since $(G_{\alpha \bar{\beta}})$ is positive definite for any $(z, v) \in T^{1,0} M \setminus \{\text{zero section}\}$, it induces a Hermitian metric $\langle \cdot, \cdot \rangle$ on the vertical subbundle $\mathcal V \subset T^{1,0} \tilde{M}$, where $\mathcal V$ is the kernel of the projection $\pi: T^{1,0} M \to M$:
\[
\langle Y, W\rangle:=G_{\alpha \overline{\beta}}(v) Y^{\alpha} \overline{W}^{\beta},\quad \forall Y, W \in \mathcal{V}.
\]

For complex Finsler manifolds, there are three equivalent notions \cite{AP}: the complex horizontal bundle, the complex horizontal map, and the complex nonlinear connection. The complex nonlinear connection is defined as follows:

\begin{definition}
A complex non-linear connection is a map
\[
\tilde{D}: \mathcal{X}\left(T^{1,0} M\right) \to \mathcal{X}\left(T_{\mathbb{C}}^{*} M \otimes T^{1,0} M\right)
\]
such that
\begin{equation}
\overline{D} \circ J \equiv 0 ; \label{2.1}
\end{equation}
\begin{equation}
\tilde{D}'' \xi \equiv 0 \quad \forall \xi \in \mathcal{X}_{\mathcal{O}}\left(T^{1,0} M\right) ;  \label{2.2}
\end{equation}
and
\[
\tilde{D} \xi_{p}'-\tilde{D} \xi_{p}=\iota_{v}^{-1} \circ\left(d^{1,0} \xi_{p}'-d^{1,0} \xi_{p}\right)
\]
for all $p \in M$ and $\xi' \in \mathcal{X}(T^{1,0} M)$ with $\xi(p)=\xi'(p)=v$.

Notice that \eqref{2.2} is equivalent to
\[
\tilde{D} \xi_{p} \circ J=J \circ \tilde{D} \xi_{p}
\]
for all $\xi \in \mathcal{X}(T^{1,0} M)$ which are holomorphic in a neighborhood of $p \in M$.
\end{definition}

A complex horizontal bundle and a complex nonlinear connection can be canonically associated to a good complex vertical connection, which is defined by \cite{AP}:

\begin{definition}
A complex vertical connection on a complex manifold is a linear connection on the complexified vertical bundle $\mathcal V_{\mathbb{C}}$, that is, a linear map $D: \Gamma(\mathcal V_{\mathbb{C}}) \to \Gamma(T_{\mathbb{C}}^{*} \tilde{M} \otimes \mathcal V_{\mathbb{C}})$ such that
\[
D(fV)=f D V+df \otimes V
\]
for every $f \in C^{\infty}(\tilde{M})$ and $V \in \Gamma(\mathcal V_{\mathbb{C}})$;
\[
D(\overline{V})=\overline{D V},\quad \forall V \in \Gamma\left(\mathcal{V}_{\mathbb{C}}\right) ;
\]
\[
D(J V)=J \circ D V,\quad \forall V \in \Gamma\left(\mathcal{V}^{1,0}\right) ;
\]
and
\[
(D V) \circ J=J \circ(D V),\quad \forall V \in \Gamma_{\mathcal{O}}\left(\mathcal{V}^{1,0}\right) .
\]

It is called \emph{good} if the bundle map $\Lambda: T^{1,0} \tilde{M} \to \mathcal V$ defined by
\[
\Lambda(X)=\nabla_{X} \iota
\]
is a bundle isomorphism.
\end{definition}

Locally, a complex vertical connection is given by
\[
D V=\left[d V^{\beta}+V^{\alpha} \omega_{\alpha}^{\beta}\right] \otimes \dot{\partial}_{\beta}
\]
for all $V \in \Gamma(\mathcal V^{1,0})$, where the $\omega_{\alpha}^{\beta}$ are local $(1,0)$-forms. In particular,
\[
\nabla_{X} V=\left[X\left(V^{\alpha}\right)+\omega_{\beta}^{\alpha}(X) V^{\beta}\right] \dot{\partial}_{\alpha},\quad \nabla_{X} \overline{V}=\overline{\nabla_{\overline{X}} V}
\]
for all $X \in T^{1,0} \tilde{M}$.

In local coordinates,
\[
\Lambda (X)=\left[ \dot{X}^{\alpha }+\omega _{\beta }^{\alpha }(X)v^{\beta }\right] \dot{\partial}_{\alpha}.
\]

The complex vertical connection $D$ is good if and only if $\mathcal H_{\mathbb{C}}=\mathcal H \oplus \overline{\mathcal H}$ is a complex horizontal bundle, where $\mathcal H=\ker \Lambda \subset T^{1,0} \tilde{M}$. In local coordinates, $D$ is good iff the matrix
\[
L_{\beta}^{\alpha}=\delta_{\beta}^{\alpha}+\tilde{\Gamma}_{\gamma \beta}^{\alpha} v^{\gamma}
\]
is invertible, where $\tilde{\Gamma}_{\gamma \beta}^{\alpha}$ are the vertical connection coefficients of $D$. If $((L^{-1})_{\nu}^{\mu})$ denotes the inverse matrix, a local frame for $H$ is given by
\[
\delta _{\alpha }=\partial _{\alpha }-\left( L^{-1}\right) _{\nu }^{\mu }\tilde{\Gamma }_{\gamma ;\alpha }^{\nu }v^{\gamma }\dot{\partial}_{\mu }=\partial _{\alpha }-\Gamma _{\alpha }^{\mu }\dot{\partial}_{\mu }
\]
where
\[
\Gamma _{\alpha }^{\mu }=\left( L^{-1}\right) _{\nu }^{\mu }\tilde{\Gamma }_{\gamma ;\alpha }^{\nu }v^{\gamma }
\]
are the Christoffel symbols of the complex non-linear connection associated to $\mathcal H$.

There are many good complex vertical connections associated to a strongly pseudoconvex complex Finsler metric on $M$, such as the Chern-Finsler connection, the complex Rund connection, the complex Berwald connection, and the complex Hashiguchi connection. We recall the first two here.

The Chern-Finsler connection is given by \cite{AP}:

\begin{theorem}
Let $F: T^{1,0} M \to \mathbb{R}^{+}$ be a strongly pseudoconvex complex Finsler metric. Then there is a unique good complex vertical connection $\tilde{D}: \Gamma(\mathcal V) \to \Gamma(T^{*} \tilde{M} \otimes \mathcal V)$ such that for all $X$ and $V, W \in \Gamma(\mathcal V)$,
\[
X\langle V,W\rangle =\left\langle \tilde{\nabla}_{X}V,W\right\rangle +\left\langle V,\tilde{\nabla}_{X}W\right\rangle .
\]
\end{theorem}

And the complex Rund connection:
\begin{theorem}[\cite{SZ}]
Let $F$ be a strongly pseudoconvex complex Finsler metric on a complex manifold $M$, and let $\langle \cdot, \cdot \rangle$ be the Hermitian structure on $\mathcal V^{1,0}$ induced by $F$. Then there is a unique complex vertical connection $\check{D}: \Gamma(\mathcal V^{1,0}) \to \Gamma(T_{\mathbb{C}}^{*} \tilde{M} \otimes \mathcal V^{1,0})$ such that:
\begin{enumerate}
\item[(i)] $\check{D}$ is good;
\item[(ii)] for any $\mathcal H \in \check{\mathcal H}^{1,0}$ and $V, W \in \mathcal V^{1,0}$,
\[
H\langle V, W\rangle=\left\langle\check{D}_{H} V, W\right\rangle+\left\langle V, \check{D}_{\overline{H}} W\right\rangle,
\]
where $\check{\mathcal H}^{1,0}$ is the $(1, 0)$ part of the complex horizontal bundle $\check{\mathcal H}_{\mathbb{C}}$ associated to $\check{D}$;
\item[(iii)] for any $V \in V^{1,0}$ and $H \in \check{\mathcal H}^{1,0}$,
\[
\dot{\theta}(V, H)=0,
\]
where $\dot{\theta}$ is the $(2, 0)$-torsion of $\check{D}$.
\end{enumerate}
\end{theorem}

In the following we denote the connection 1-forms with respect to the Cartan connection, the Chern-Finsler connection and the complex Rund connection by $\hat{\omega}_{j}^{i}$, $\tilde{\omega}_{j}^{i}$, $\check{\omega}_{j}^{i}$ respectively. Locally they are expressed as
\[
\hat{\omega}_{j}^{i}=\hat{\Gamma}_{j;k}^{i}dz^{k}+\hat{\Gamma}_{jk}^{i}\hat{\psi}^{k},\quad \tilde{\omega}_{j}^{i}=\tilde{\Gamma}_{j;k}^{i}dz^{k}+\tilde{\Gamma}_{jk}^{i}\psi^{k},\quad \check{\omega}_{j}^{i}=\check{\Gamma}_{j;k}^{i}dz^{k}.
\]

Here $\hat{\Gamma}_{jk}^{i}$, $\tilde{\Gamma}_{j;k}^{i}$, and $\check{\Gamma}_{j;k}^{i}$ are called the horizontal connection coefficients of the corresponding connections. It is known that $\check{\Gamma}_{j;k}^{i}=\tilde{\Gamma}_{j;k}^{i}$, i.e., the complex Rund connection is just the horizontal part of the Chern-Finsler connection.

The Hermitian connection is parallel with respect to the corresponding Hermitian metric, and it is torsion-free when the Hermitian metric is a K\"ahler metric. However, in the complex Finsler case, there exists no complex vertical connection satisfying both the metric compatibility condition and the torsion-free condition at the same time, even in the K\"ahler-Finsler case, and the torsion is usually more complicated. For example, the torsion of the Chern-Finsler connection is given by \cite{AP}:
\[
T=\theta^{\sigma} \otimes \delta_{\sigma}+\tau^{\alpha} \dot{\partial}_{\alpha},
\]
\[
\theta^{\sigma}=\frac{1}{2}\left[\Gamma_{\nu ; \mu}^{\sigma}-\Gamma_{\mu ; \nu}^{\sigma}\right] d z^{\mu} \wedge d z^{\nu}+\Gamma_{\nu \gamma}^{\sigma} \psi^{\gamma} \wedge d z^{\nu},
\]
\[
\tau ^{\alpha }=-\delta _{\overline{\nu }}\left( \Gamma _{;\mu }^{\alpha }\right) dz^{\mu }\wedge d\overline{z}^{\nu }-\Gamma _{\overline{\beta};\mu }^{\alpha }dz^{\mu }\wedge \overline{\psi ^{\beta }}.
\]

More generally, the torsion $T=\theta+\tau$ of a good complex vertical connection is given by [10]:
\[
\begin{aligned}
\theta=&\left[\frac{1}{2}\left(\Gamma_{\nu ; \sigma}^{\mu}-\Gamma_{\sigma ; \nu}^{\mu}\right) d z^{\sigma} \wedge d z^{\nu}+\Gamma_{\nu \gamma}^{\mu} \psi^{\gamma} \wedge d z^{\nu}\right] \otimes \delta_{\mu}+\left\{\frac{1}{2}\left[\delta_{\mu}\left(\Gamma_{\nu}^{\alpha}\right)-\delta_{\nu}\left(\Gamma_{\mu}^{\alpha}\right)\right] d z^{\mu} \wedge d z^{\nu}\right. \\
&\left.+\left[\dot{\partial}_{\beta}\left(\Gamma_{\mu}^{\alpha}\right)-\Gamma_{\beta ; \mu}^{\alpha}\right] \psi^{\beta} \wedge d z^{\mu}+\frac{1}{2}\left(\Gamma_{\beta \gamma}^{\alpha}-\Gamma_{\gamma \beta}^{\alpha}\right) \psi^{\beta} \wedge \psi^{\gamma}\right\} \otimes \dot{\partial}_{\alpha}, \\
\tau=&-\left[\delta_{\overline{\nu}}\left(\Gamma_{\mu}^{\alpha}\right) d z^{\mu} \wedge d \overline{z}^{\nu}+\dot{\partial}_{\overline{\beta}}\left(\Gamma_{\mu}^{\alpha}\right) d z^{\mu} \wedge \overline{\psi^{\beta}}\right] \otimes \dot{\partial}_{\alpha}.
\end{aligned}
\]

By \cite{AP}, there are three notions of K\"ahlerian Finsler metrics:

\begin{definition}
A strongly pseudoconvex Finsler metric $F$ is called \emph{strongly K\"ahler} if the horizontal part of the $(2, 0)$-torsion vanishes, that is, iff
\[
\forall H, K \in \mathcal{H},\quad \tilde{\theta}(H, K)=0 ;
\]
\emph{K\"ahler} if
\[
\forall H \in \tilde{\mathcal{H}},\quad \tilde{\theta}(H, \chi)=0,
\]
and \emph{weakly K\"ahler} if
\[
\forall H \in \tilde{\mathcal{H}},\quad \langle\tilde{\theta}(H, \chi), \chi\rangle=0 .
\]
\end{definition}

By \cite{CS}, the definition of strongly K\"ahler Finsler metric is equivalent to the K\"ahler Finsler metric. In local coordinates, $F$ is strongly K\"ahler iff
\[
\Gamma_{\mu ; \nu}^{\alpha}=\Gamma_{\nu ; \mu}^{\alpha} ;
\]
it is K\"ahler iff
\[
\tilde{\Gamma}_{\mu ; \nu}^{\alpha} v^{\mu}=\tilde{\Gamma}_{\nu ; \mu}^{\alpha} v^{\mu} ;
\]
it is weakly K\"ahler iff
\[
G_{\alpha}\left[\tilde{\Gamma}_{\mu ; \nu}^{\alpha}-\tilde{\Gamma}_{\nu ; \mu}^{\alpha}\right] v^{\mu}=0,
\]
where $\chi:=v^i\delta_i$ is the horizontal radial vector field.

$F$ is called a \emph{complex Berwald metric} if $\tilde{\Gamma}_{\nu ; \mu}^{\alpha}$ are independent of the tangent direction $v=(v_{1}, \cdots, v_{n})$; $F$ is called a \emph{K\"ahler-Berwald metric} if $F$ is both a K\"ahler-Finsler metric and a complex Berwald metric. $F$ is called a \emph{real Berwald metric} if $\hat{\Gamma}_{i ; j}^{k}$ are independent of the tangent direction.

%\section{Equivalent conditions of K\"ahler Finsler metric}\label{sec3}
%In this section, we study equivalent conditions for a complex Finsler metric to be a K\"ahler Finsler metric.
%In the strongly convex real Berwald case, we prove that if the complex structure is parallel with respect to the nonlinear connection of the Cartan connection, then the metric is necessarily a K\"ahler Berwald metric.
%Next, we consider cases without the real Berwald condition;
%Finally, we consider the strongly pseudoconvex case.

\section{Real Berwald strongly convex complex Finsler metrics with parallel complex structure}\label{sec3}

Now let $(M, J, F, \hat{\nabla})$ be a strongly convex complex Finsler manifold, where $F$ is a strongly convex complex Finsler, real Berwald metric. Here $\hat{\nabla}$ is the nonlinear connection corresponding to the Cartan connection of $F$, and the complex structure $J$ of $M$ is parallel with respect to the nonlinear connection $\hat{\nabla}$, i.e. $\hat{\nabla} J=0$.

Since $F$ is real Berwald, in local coordinates we have
\begin{equation}
\hat{\nabla} X=\left(d X^{i}+X^{j} \hat{\Gamma}_{j ; k}^{i}(x) d x^{k}\right) \otimes \partial_{i}^{o},\quad \forall X \in \Gamma\left(\tilde{M}_{\mathbb{R}}\right), \label{3.1}
\end{equation}
where $\hat{\Gamma}_{j ; k}^{i}(x)=\hat{\Gamma}_{k ; j}^{i}(x)$ depends only on $x$. The condition $\hat{\nabla} J=0$ is equivalent to
\[
\hat{\nabla}(J X)=J(\hat{\nabla} X),\quad \forall X \in \Gamma\left(\tilde{M}_{\mathbb{R}}\right).
\]

Thus $\hat{\nabla}_{\partial_{j}^{o}}(J \partial_{\alpha}^{o})=J(\hat{\nabla}_{\partial_{j}^{o}} \partial_{\alpha}^{o})$. By
\[
\hat{\nabla}_{\partial_{j}^{o}}\left(J \partial_{\alpha}^{o}\right)=\hat{\nabla}_{\partial_{j}^{o}}\left(\partial_{\alpha+n}^{o}\right)=\hat{\Gamma}_{\alpha+n ; j}^{i} \partial_{i}^{o},
\]
\[
\begin{aligned}
J\left( \hat{\nabla}_{\partial_{j}^{o}}\partial_{\alpha }^{o}\right) &=J\left( \hat{\Gamma }_{\alpha ;j}^{i}\partial_{i}^{o}\right) =\hat{\Gamma }_{\alpha ;j}^{\beta }\partial_{\beta +n}^{o}-\hat{\Gamma }_{\alpha ;j}^{\beta +n}\partial_{\beta }^{o},
\end{aligned}
\]
we obtain
\[
\hat{\Gamma }_{\alpha ;j}^{\beta }=\hat{\Gamma }_{\alpha +n;j}^{\beta +n};\quad \hat{\Gamma }_{\alpha +n;j}^{\beta }=-\hat{\Gamma }_{\alpha ;j}^{\beta +n}.
\]

Define
\begin{equation}
\Gamma_{\alpha ; \gamma}^{\beta}=\hat{\Gamma}_{\alpha ; \gamma}^{\beta}+i \hat{\Gamma}_{\alpha ; \gamma}^{\beta+n}=\hat{\Gamma}_{\alpha ; \gamma+n}^{\beta+n}-i \hat{\Gamma}_{\alpha ; \gamma+n}^{\beta}=-i \Gamma_{\alpha ; \gamma+n}^{\beta}. \label{3.6}
\end{equation}
Then for any $\eta=\eta^{\gamma} \partial_{\gamma} \in \Gamma(\tilde{M})$, and $\xi=\xi^{c} \partial_{c}^{o}=\eta+\bar{\eta} \in \Gamma(\tilde{M}_{\mathbb{R}})$, we have
\[
\xi^{c} \Gamma_{\alpha ; c}^{\beta}=\xi^{\gamma} \hat{\Gamma}_{\alpha ; \gamma}^{\beta}+\xi^{\gamma+n} \hat{\Gamma}_{\alpha ; \gamma+n}^{\beta}=\left(\xi^{\gamma}+i \xi^{\gamma+n}\right) \hat{\Gamma}_{\alpha ; \gamma}^{\beta}=\eta^{\gamma} \Gamma_{\alpha ; \gamma}^{\beta}.
\]

We also define
\begin{equation}
\Gamma _{\alpha }^{\beta }(z):=\hat{\Gamma }_{\alpha }^{\beta }+i\hat{\Gamma }_{\alpha }^{\beta +n}=u^{c}\left( \hat{\Gamma }_{\alpha ;c}^{\beta }+i\hat{\Gamma }_{\alpha ;c}^{\beta +n}\right) =u^{c}\Gamma _{\alpha ;c}^{\beta }=v^{\gamma }\Gamma _{\alpha ;\gamma }^{\beta },\quad v\in T_{z}^{1,0}M. \label{3.8}
\end{equation}

We can now define a complex vertical connection:

\begin{theorem}
Let $(M, J, F, \hat{\nabla})$ be a strongly convex real Berwald complex Finsler manifold, with the complex structure $J$ parallel with respect to the nonlinear connection $\hat{\nabla}$. Define $D: \Gamma(\mathcal V) \to \Gamma(T^{*} \tilde{M} \otimes \mathcal V)$ by
\begin{equation}
D(Y)=\left(d Y^{\alpha}+Y^{\beta} \Gamma_{\beta ; \gamma}^{\alpha} d z^{\gamma}\right) \otimes \dot{\partial}_{\alpha},\quad \forall Y=Y^{\alpha} \dot{\partial}_{\alpha} \in \Gamma(\mathcal{V}), \label{connection}
\end{equation}
where $\Gamma_{\beta ; \gamma}^{\alpha}$ are given by \eqref{3.6}. Then $D$ is a good complex vertical connection.
\end{theorem}

\begin{proof}
By definition, $D$ is clearly a complex vertical connection. Notice that the vertical connection coefficients $\Gamma_{\beta \gamma}^{\alpha}$ vanish, thus the matrix $(L_{\beta}^{\alpha})=(\delta_{\beta}^{\alpha}+\Gamma_{\gamma \beta}^{\alpha} v^{\gamma})=I$ is clearly invertible. Hence $D$ is a good complex vertical connection.
\end{proof}

A local frame for the horizontal bundle $H$ associated to $D$ is given by
\[
\delta _{\alpha }=\partial _{\alpha }-\left( L^{-1}\right) _{\nu }^{\mu }\Gamma _{\gamma ;\alpha }^{\nu }v^{\gamma }\dot{\partial}_{\mu }=\partial _{\alpha }-\Gamma _{\alpha }^{\mu }\dot{\partial}_{\mu }
\]
where $\Gamma_{\alpha}^{\mu}$ is given by \eqref{3.8} as the Christoffel symbols of the complex non-linear connection associated to $\mathcal H$.

Indeed, the complex horizontal bundle $\mathcal H$ is exactly the $(1,0)$ part of the complexified bundle $(\mathcal H_{\mathbb{R}})_{\mathbb{C}}$ of the real horizontal bundle $\mathcal H_{\mathbb{R}}$, since we have:
\[
2 \delta _{\alpha }=\partial _{\alpha }^{o}-i\partial _{\alpha +n}^{o}-(\hat{\Gamma }_{\alpha }^{\mu }+i\hat{\Gamma }_{\alpha }^{\mu +n})(\dot{\partial}_{\mu }^{o}-i \dot{\partial}_{\mu +n}^{o})=\partial _{\alpha }^{o}-\hat{\Gamma }_{\alpha }^{b} \dot{\partial}_{b}^{o}-i\left( \partial _{\alpha +n}^{o}-\hat{\Gamma }_{\alpha +n}^{b} \dot{\partial}_{b}^{o}\right) =\delta _{\alpha }^{o}-i \delta _{\alpha +n}^{o}.
\]

We then have the following main theorem:

\begin{theorem}
Let $(M, J, F, \hat{\nabla})$ be a strongly convex real Berwald complex Finsler manifold, with the complex structure $J$ parallel with respect to the nonlinear connection $\hat{\nabla}$. The good complex vertical connection $D$ defined by \eqref{connection} is the complex Rund connection, and $F$ is a K\"ahler Berwald manifold.
\end{theorem}

\begin{proof}
The fundamental tensor is given by
\[
\tilde{\gamma}_{\alpha \overline{\beta}}=\frac{1}{4}\left(\dot{\partial}_{\alpha}^{o}-i \dot{\partial}_{\alpha+n}^{o}\right)\left(\dot{\partial}_{\beta}^{o}+i \dot{\partial}_{\beta+n}^{o}\right) G=G_{u^{\alpha} u^{\beta}}+G_{u^{\alpha+n} u^{\beta+n}}+i\left(G_{u^{\alpha} u^{\beta+n}}-G_{u^{\alpha+n} u^{\beta}}\right),
\]
where $G_{u^{a} u^{b}}:=\frac{\partial^{2} G}{\partial u^{a} \partial u^{b}}$. Thus the Hermitian metric on $V$ and the Riemannian metric on $\mathcal{V}_{\mathbb{R}}$ satisfy the relation:
\[
\begin{aligned}
4\langle Y,W\rangle&=4 G_{\alpha \overline{\beta}} Y^{\alpha} \overline{W}^{\beta} \\
&=\left[ G_{u^{\alpha}u^{\beta}}+G_{u^{\alpha+n}u^{\beta+n}}+i\left( G_{u^{\alpha}u^{\beta+n}}-G_{u^{\alpha+n}u^{\beta}}\right) \right] (Y^{\alpha}+iY^{\alpha+n})(\overline{W}^{\beta}-i\overline{W}^{\beta+n}) \\
&=\left[ G_{u^{\alpha}u^{\beta}}+G_{u^{\alpha+n}u^{\beta+n}}\right] \left( Y^{\alpha}W^{\beta}+Y^{\alpha+n}W^{\beta+n}\right) +\left( G_{u^{\alpha}u^{\beta+n}}-G_{u^{\alpha+n}u^{\beta}}\right) \left( Y^{\alpha+n}W^{\beta}-Y^{\alpha}W^{\beta+n}\right) \\
&\quad +i\left[ \left( G_{u^{\alpha}u^{\beta+n}}-G_{u^{\alpha+n}u^{\beta}}\right) \left( Y^{\alpha}W^{\beta}+Y^{\alpha+n}W^{\beta+n}\right) -\left( G_{u^{\alpha}u^{\beta}}+G_{u^{\alpha+n}u^{\beta+n}}\right) \left( Y^{\alpha+n}W^{\beta}-Y^{\alpha}W^{\beta+n}\right) \right].
\end{aligned}
\]
Here we denote $\langle A+iB | C+iD\rangle:=\langle A | C\rangle+i\langle A | D\rangle+i\langle B | C\rangle-\langle B | D\rangle$ by complex linear extension. Thus
\begin{equation}
\delta_{\eta}\left\langle Y, W\right\rangle=2 \delta_{\eta}\left\langle Y | \overline{W}\right\rangle=2\left\langle\hat{\nabla}_{\delta_{\eta}} Y | \overline{W}\right\rangle+2\left\langle Y | \hat{\nabla}_{\delta_{\eta}} \overline{W}\right\rangle, \label{3.13}
\end{equation}
where we define $\hat{\nabla}_{A+i B}(C+i D):=\hat{\nabla}_{A}(C)+i \hat{\nabla}_{A}(D)+i \hat{\nabla}_{B}(C)-\hat{\nabla}_{B}(D)$.

We have
\begin{equation}
\hat{\nabla}_{\delta_{\eta}} Y=\delta_{\eta}\left(Y^{\alpha}\right) \dot{\partial}_{\alpha}+Y^{\alpha} \hat{\nabla}_{\delta_{\eta}} \dot{\partial}_{\alpha} . \label{3.14}
\end{equation}

A direct computation shows:
\begin{equation}
\begin{aligned}
4\hat{\nabla}_{\delta_{\eta}} \dot{\partial}_{\alpha }&=\hat{\nabla}_{\delta_{\eta }^{o}-i\delta_{\eta +n}^{o}}(\dot{\partial}_{\alpha }^{o}-i\dot{\partial}_{\alpha +n}^{o})=\left[ \left( \hat{\Gamma }_{\alpha ;\eta }^{b}-\hat{\Gamma }_{\alpha +n;\eta +n}^{b}\right) -i\left( \hat{\Gamma }_{\alpha +n;\eta }^{b}+\hat{\Gamma }_{\alpha ;\eta +n}^{b}\right) \right] \dot{\partial}_{b}^{o} \\
&=\left( \hat{\Gamma }_{\alpha ;\eta }^{\beta }-\hat{\Gamma }_{\alpha +n;\eta +n}^{\beta }\right) \dot{\partial}_{\beta }^{o}+\left( \hat{\Gamma }_{\alpha ;\eta }^{\beta +n}-\hat{\Gamma }_{\alpha +n;\eta +n}^{\beta +n}\right) \dot{\partial}_{\beta +n}^{o} \\
&\quad -i\left( \hat{\Gamma }_{\alpha +n;\eta }^{\beta }+\hat{\Gamma }_{\alpha ;\eta +n}^{\beta }\right) \dot{\partial}_{\beta }^{o}-i\left( \hat{\Gamma }_{\alpha +n;\eta }^{\beta +n}+\hat{\Gamma }_{\alpha ;\eta +n}^{\beta +n}\right) \dot{\partial}_{\beta +n}^{o} \\
&=2\Gamma _{\alpha ;\eta }^{\beta }\dot{\partial}_{\beta }^{o}-2i\Gamma _{\alpha ;\eta }^{\beta }\dot{\partial}_{\beta +n}^{o}=4\Gamma _{\alpha ;\eta }^{\beta }\dot{\partial}_{\beta }.
\end{aligned} \label{3.15}
\end{equation}

Combining \eqref{3.14} and \eqref{3.15}, we get
\[
\nabla_{\delta_{\eta}} Y=\delta_{\eta}\left(Y^{\alpha}\right) \dot{\partial}_{\alpha}+Y^{\alpha} \nabla_{\delta_{\eta}} \dot{\partial}_{\beta}=\hat{\nabla}_{\delta_{\eta}} Y.
\]

From \eqref{3.13} it follows that
\[
\delta_{\eta}\left\langle Y, W\right\rangle=2\left\langle\nabla_{\delta_{\eta}} Y | \overline{W}\right\rangle+2\left\langle Y | \nabla_{\delta_{\eta}} \overline{W}\right\rangle=\left\langle\nabla_{\delta_{\eta}} Y, W\right\rangle+\left\langle Y | \nabla_{\overline{\delta}_{\eta}} W\right\rangle.
\]

Since $D$ is a good complex vertical connection with $\Gamma_{\alpha ; \eta}^{\beta}=\Gamma_{\eta ; \alpha}^{\beta}$, $\Gamma_{\eta}^{\beta}(v)=\Gamma_{\eta ; \alpha}^{\beta}(z) v^{\eta}$ for all $v \in T^{1,0} M$, and $\Gamma_{\alpha \eta}^{\beta}=0$, by the torsion formula the $(2,0)$-torsion of $D$ is given by
\[
\frac{1}{2}\left[\delta_{\mu}\left(\Gamma_{\nu}^{\alpha}\right)-\delta_{\nu}\left(\Gamma_{\mu}^{\alpha}\right)\right] d z^{\mu} \wedge d z^{\nu} \otimes \dot{\partial}_{\alpha}.
\]
Hence for any $V \in \mathcal V^{1,0}$ and $H \in \mathcal H^{1,0}$, we have $\dot{\theta}(V, H)=0$.

By Theorem 2.3, the connection $D$ is exactly the complex Rund connection. Moreover, $\Gamma_{\alpha ; \eta}^{\beta}$ coincide with the horizontal Chern-Finsler connection coefficients, thus the metric $F$ is a K\"ahler Berwald metric.
\end{proof}

Below we give some examples satisfying the conditions in Theorem 3.2.

By Theorem 4.3 in \cite{DH}, on an (compact or noncompact) irreducible Hermitian symmetric space $M=G / H$, for any $G$-invariant Finsler metric on $M$, there exists an invariant complex structure $J$ on $M$ such that $(M, J, F)$ is a symmetric complex Finsler space. By construction, $J$ is $G$-invariant and $F$ is $J$-invariant.

For any $X \in T_{o} M$, there is a geodesic curve given by the exponential map:
\[
\gamma_{X}(t):=\operatorname{Exp}_{o}(t X),\quad t \in \mathbb{R}.
\]

Since at any point $x \in M$ there exists a holomorphic involutive isometry $s_{x}$, we can define $T_{t}:=s_{\gamma(\frac{t}{2})} s_{o}$, a one-parameter subgroup in $G$. Since $s_x$ preserves the complex structure, so does $d T_{t}$. Thus for the nonlinear connection $\hat{\nabla}$ associated to the Cartan connection and any $Y \in T_o M$, we have:
\begin{equation}
\begin{aligned}
\hat{\nabla}_{\dot{\gamma}_{X}(t)}\left[ (dT_{t})_{o}(J_{o}Y)\right] &=\hat{\nabla}_{\dot{\gamma}_{X}(t)}\left[ J_{\gamma_{X}(t)}\circ (dT_{t})_{o}(Y)\right] \\
&=\left( \hat{\nabla}_{\dot{\gamma}_{X}(t)}J_{\gamma_{X}(t)}\right) \left[ (dT_{t})_{o}(Y)\right] +J_{\gamma_{X}(t)}\circ \hat{\nabla}_{\dot{\gamma}_{X}(t)}\left[ (dT_{t})_{o}(Y)\right] .
\end{aligned} \label{3.19}
\end{equation}

Since $(d T_{t})_{o}$ is parallel along $\gamma_X$, we have $\hat{\nabla}_{\dot{\gamma}_{X}(t)}[(d T_{t})_{o}(J_{o} Y)]=\hat{\nabla}_{\dot{\gamma}_{X}(t)}[(d T_{t})_{o}(Y)]=0$. From \eqref{3.19} we deduce $(\hat{\nabla}_{\dot{\gamma}_{X}(t)} J_{\gamma_{X}(t)})[(d T_{t})_{o}(Y)]=0$. By the arbitrariness of $X, Y \in \Gamma(T_{o} M)$, we have $\hat{\nabla} J=0$. Thus $(M, J, F)$ satisfies the condition of Theorem 3.2, and we obtain:

\begin{corollary}
On an (compact or noncompact) irreducible Hermitian symmetric space $(M=G / H, J, F)$, any $G$-invariant strongly convex complex Finsler metric $F$ is a K\"ahler Berwald metric.
\end{corollary}

Finally, by considering the equivalent condition for the parallelism of the complex structure $J$ with respect to the nonlinear connection $\hat{\nabla}$, we can obtain further equivalent conditions for a strongly convex real Berwald metric to be K\"ahler Berwald.

\begin{remark}
By the proof given in this section, we indeed proved that $\hat\nabla J=0$ implies $F$ being K\"ahler--Finsler without imposing the real Berwald condition. However, for the purpose of preserving the form of the original 2022 manuscript, no modifications have been made to the proof, and the statement about results in this section.
\end{remark}

\section{Holomorphically invariant complex Finsler metrics of $\mathfrak R_{IV}$}\label{4}

	\subsection{Bergmann metric on the classic domain of type IV: $\mathfrak{R}_{IV}$}
	
	The classical domain of type $\mathfrak{R}_{IV}$ consisting points $z=(z^1,\cdots,z^n)\in \mathbb{C}^n$ satisfying
	\begin{equation}
		\mathfrak{R}_{IV}: \left\{
		\begin{array}{ll}
			|zz'|<1, &\\
			1+|zz'|^2-2z\overline{z'}>0, &
		\end{array}
		\right. \quad n\geq 5.
	\end{equation}
	
	For each  $z_0\in\mathfrak{R}_{IV}$, there exists $\phi_{z_0}\in\mbox{Aut}(\mathfrak{R}_{IV})$
	such that $\phi_{z_0}(z_0)=0$ and for $z\in\mathfrak{R}_{IV}$ \cite{Look}
	\begin{equation}
		\phi_{z_0}(z)=\left\{\left[\left(\frac{1+zz'}{2},\frac{1-zz'}{2i}\right)-zX_0'\right]A\begin{pmatrix}
			1 \\
			i \\
		\end{pmatrix}
		\right\}^{-1}\left[z-\left(\frac{1+zz'}{2},\frac{1-zz'}{2i}\right)X_0\right]D,
		\label{phi-z0}
	\end{equation}
	where
	\begin{eqnarray}
		X_0&=&2\begin{pmatrix}
			z_0z_0'+1 & i(z_0z_0'-1) \\
			\overline{z_0z_0'}+1 & -i(\overline{z_0z_0'}-1) \\
		\end{pmatrix}^{-1}
		\begin{pmatrix}
			z_0 \\
			\overline{z_0} \\
		\end{pmatrix}\nonumber\\
		&=&\frac{-1}{1-|z_0z_0'|^2}\begin{pmatrix}
			(\overline{z_0z_0'}-1)z_0+(z_0z_0'-1)\overline{z_0}  \\
			i(z_0z_0'+1)\overline{z_0}-i(\overline{z_0z_0'}+1)z_0 \\
		\end{pmatrix}
		\label{X0}
	\end{eqnarray}
	is a real $2$-by-$n$ matrix satisfying $I-X_0X_0'>0$,
	$A$ is a $2$-by-$2$ real matrix and $D$ is an $n$-by-$n$ real matrix satisfying
	\begin{equation}
		AA'=(I-X_0X_0')^{-1},\quad DD'=(I-X_0'X_0)^{-1},\quad \det A>0.\label{D}
	\end{equation}
	The isotropic subgroup $\mbox{Iso}(\mathfrak{R}_{IV})$ of $\mbox{Aut}(\mathfrak{R}_{IV})$ at the origin $0\in \mathfrak{R}_{IV}$ is given by
	\begin{equation}
		w=e^{i\theta}z D,\quad DD'=I.
	\end{equation}
	
	By \eqref{X0}, we have
	\begin{eqnarray}
		I_n-X_0'X_0=I_n+c_1z_0'z_0-c_2z_0'\overline{z_0}-c_2\overline{z_0'}z_0+\overline{c_1}\overline{z_0'z_0},\label{IXX0}
	\end{eqnarray}
	where
	\begin{equation}
		c_1=\frac{4\overline{z_0z_0'}}{(1-|z_0z_0'|^2)^2},\quad c_2=\frac{2(1+|z_0z_0'|^2)}{(1-|z_0z_0'|^2)^2}.
	\end{equation}
	Substituting \eqref{IXX0} into the second equality of \eqref{D} and by a linear algebra argument, we obtain
	\begin{equation}
		DD'=I_n+\frac{2}{1+|z_0z_0'|^2-2\|z_0\|^2}(z_0'\overline{z_0}+\overline{z_0'}z_0).
	\end{equation}

	Denoting by \begin{equation} \Delta(z)=1+|zz'|^2-2|z|^2, \label{delta}\end{equation}
	then
	\begin{equation}
		DD'=I_n+\frac{2}{\Delta(z_0)}(z_0'\overline{z_0}+\overline{z_0'}z_0).
	\end{equation}

	Taking
	\begin{equation}
		A=\frac{1}{2\Delta(z_0)^{\frac{1}{2}}}\begin{pmatrix}
			-i(z_0z_0'-\overline{z_0z_0'}) & z_0z_0'+\overline{z_0z_0'}-2 \\
			z_0z_0'+\overline{z_0z_0'}+2& i(z_0z_0'-\overline{z_0z_0'}) \\
		\end{pmatrix},
		\label{A}
	\end{equation}
	and
	\begin{eqnarray*}
		D=D'&=&\Bigg[I_n+\frac{4}{\Delta(z_0)}\begin{pmatrix}
			\mbox{Re}(z_0) \\
			\mbox{Im}(z_0)\\
		\end{pmatrix}'
		\begin{pmatrix}
			\mbox{Re}(z_0) \\
			\mbox{Im}(z_0) \\
		\end{pmatrix}
		\Bigg]^{\frac{1}{2}}>0,
	\end{eqnarray*}
	then $A$ and $D$ satisfy \eqref{D}.
	
	Differentiating \eqref{phi-z0} and then setting $z=z_0$, one has
	\begin{equation}
		d\phi_{z_0}(z)|_{z=z_0}=\left\{dzD-dzz_0'(1,i)X_0D\right\}\left\{\left[\left(\frac{1+z_0z_0'}{2},\frac{1-z_0z_0'}{2i}\right)-z_0X_0'\right]A\begin{pmatrix}
			1 \\
			i \\
		\end{pmatrix}
		\right\}^{-1}.
	\end{equation}
	Substituting \eqref{A} and \eqref{X0} into the above equation, one gets \cite{Look}
	\begin{equation}
		d\phi_{z_0}(z)|_{z=z_0}=\frac{1}{i\Delta(z_0)^{\frac{1}{2}}}dz\left\{I-2\frac{z_0'\overline{z_0}-\overline{z_0z_0'}z_0'z_0}{1-|z_0z_0'|^2}\right\}D.
	\end{equation}
	Thus for any tangent vector $\pmb{v}=v^i\frac{\partial}{\partial z^i}\in T_{z_0}^{1,0}\mathfrak{R}_{IV}$, we have
	\begin{equation}
		(\phi_{z_0})_\ast (v)=\frac{1}{i\Delta(z_0)^{\frac{1}{2}}}v\left\{I-2\frac{z_0'\overline{z_0}-\overline{z_0z_0'}z_0'z_0}{1-|z_0z_0'|^2}\right\}D
	\end{equation}
	is a tangent vector at the origin since $\phi_{z_0}(z_0)=0$.
	
	Denote $f^2(v)=\mbox{tr}\{v\overline{v'}\}$ the canonical complex Eculidean norm on $\mathbb{C}^n\cong T_0^{1,0}\mathfrak{R}_{IV}$, which is both $U(n;\mathbb{C})$-invariant and $O(n;\mathbb{R})$-invariant. Then the Bergmann metric $F^2(z_0; v)$ satisfies
	\begin{equation}
		F^2(z_0;v)=F^2(\phi_{z_0}(z_0); (\phi_{z_0})_\ast(v))=F^2(0;(\phi_{z_0})_\ast(v))=f^2\left((\phi_{z_0})_\ast(v)\right).\label{F2z0}
	\end{equation}
	Note that by a direct calculation, we have
	\begin{eqnarray*}
		&&\left\{I_n-\frac{2}{1-|z_0z_0'|^2}\left[z_0'\overline{z_0}-(\overline{z_0z_0'})z_0'z_0\right]\right\}
		\left[I_n+\frac{2}{\Delta(z_0)}(z_0'\overline{z_0}+\overline{z_0'}z_0)\right]\\
		&=&I_n+\frac{2}{\Delta(z_0)}\left[\overline{z_0'}z_0-(\overline{z_0z_0'})z_0'z_0\right]
	\end{eqnarray*}
	and
	\begin{eqnarray*}
		&&\left\{I_n+\frac{2}{\Delta(z_0)}\left[\overline{z_0'}z_0-(\overline{z_0z_0'})z_0'z_0\right]\right\}\left\{I_n-\frac{2}{1-|z_0z_0'|^2}\left[z_0'\overline{z_0}-(z_0z_0')\overline{z_0'z_0}\right]\right\}\\
		&=&I_n+\frac{2}{\Delta(z_0)}\left\{
		2\|z_0\|^2z_0'\overline{z_0}
		+\left(\overline{z_0'}z_0-z_0'\overline{z_0}\right)-\left[(\overline{z_0z_0'})z_0'z_0
		+(z_0z_0')\overline{z_0'z_0}\right]
		\right\}.
	\end{eqnarray*}
	Denoting by
	\begin{equation}
		A(z_0)=\left(a_{i\overline{j}}(z_0)\right)=\frac{1}{\Delta(z_0)}\Bigg\{I_n+\frac{2}{\Delta(z_0)}\left\{
		2\|z_0\|^2z_0'\overline{z_0}
		+(\overline{z_0'}z_0-z_0'\overline{z_0})-\left[(\overline{z_0z_0'})z_0'z_0
		+(z_0z_0')\overline{z_0'z_0}\right]
		\right\}\Bigg\},\label{1.14}
	\end{equation}
	then
	$$
	F^2(z_0;v)=vA(z_0)\overline{v'}.
	$$

	Since $z_0$ is an arbitrary fixed point in $\mathfrak{R}_{IV}(n)$, changing $z_0$ to $z$, we obtain the Bergmann
	metric $ds_B^2=F_B^2(z;dz)$, where
	\begin{eqnarray*}
		F_B^2(z;v)&=&\frac{v\overline{v'}}{\Delta(z)}+\frac{2v}{\Delta(z)^2}\Big\{2|z|^2z'\overline{z}
		+(\overline{z'}z-z'\overline{z})-[(\overline{zz'})z'z
		+(zz')\overline{z'z}]
		\Big\}\overline{v'}.
	\end{eqnarray*}
	
	We also notice  that
	\begin{eqnarray}
		\left\{I_n+\frac{2}{\Delta(z_0)}\left[\overline{z_0'}z_0-(\overline{z_0z_0'})z_0'z_0\right]\right\}\left\{I_n-\frac{2}{1-|z_0z_0'|^2}\left[\overline{z_0'}z_0-(\overline{z_0z_0'})z_0'z_0\right]\right\}
		=I_n.\label{I-^{-1}}
	\end{eqnarray}
	Thus
	$$
	\left[(\phi_{z_0})_\ast(v)\right]\left[(\phi_{z_0})_\ast(v)\right]'=\frac{v v'}{\Delta(z_0)}.
	$$
	Thus we have
	\begin{proposition} The Kobayashi metric $F_K^2(z;v)$ on $\mathfrak{R}_{IV}(n)$ is given by
		\begin{eqnarray*}
			F_K^2(z;v)&=&F_B^2(z;v)+\sqrt{F_B^4(z;v)-\frac{|vv'|^2}{\Delta(z)^2}}.
		\end{eqnarray*}
	\end{proposition}

	\subsection{Proof of Theorem \ref{R4}}
	
	In this section, we consider $\mbox{Iso}_0(\mathfrak{R}_{IV})$-invariant strongly pseudoconvex and strongly convex complex Minkowski norms on  $T_0^{1,0}\mathfrak{R}_{IV}.$
	
The isotropic subgroup $\mbox{Iso}_0(\mathfrak{R}_{IV})$ of $\mbox{Aut}(\mathfrak{R}_{IV})$ at the origin $0\in \mathfrak{R}_{IV}$ is given by (cf. \cite{Look})
\begin{equation}
	w=e^{i\theta}z D,\quad DD'=I.\label{Iso}
\end{equation}
Thus an $\mbox{Iso}_0(\mathfrak{R}_{IV})$-invariant strongly pseudoconvex complex Minkowski norms $f$ on  $T_0^{1,0}\mathfrak{R}_{IV}$ is a strongly pseudoconvex complex Minkowski norm satisfy
\begin{eqnarray}
	f(v)&=&f(vD)\quad for\quad DD'=I;\label{gA}\\
	f(v)&=&f(e^{i\theta}v) \quad for\quad  \forall\theta\in \mathbb{R};\label{gB}.
\end{eqnarray}
For any nonzero tangent vector $0\neq v\in \mathbb{C}^n \cong T_{0}^{1,0} \mathfrak{R}_{IV},$ $ v $ can be writen as follows:
\begin{equation}
	v=(R+iS)x+iTy \quad with \quad xy'=0, xx'=yy'=1,\quad and \quad R>0,\,T\geq0, S\in\mathbb{R},\label{v-rst}
\end{equation}
where
\begin{eqnarray}
	R^2&=&\frac{1}{4}(v+\bar{v})(v+\bar{v})'
	=\frac{1}{4}(vv'+\overline{vv'}+2v\bar{v}'), \label{v-r} \\
	RS&=&\frac{1}{4i}(v-\bar{v})(v+\bar{v})'
	=\frac{1}{4i}(vv'-\overline{vv'}), \label{v-rs} \\
	S^2+T^2&=&-\frac{1}{4}(v-\bar{v})(v-\bar{v})'
	=\frac{1}{4}(2v\bar{v}'-vv'-\overline{vv'}). \label{v-st}
\end{eqnarray}
Take an orthogonal matrix $ O(v) $ by using
$x,y $ %(resp. $ (y,x) $)
as its first two columns, thus we have
\begin{equation} vO(v)=\big(R+iS,iT,0,\cdots,0\big)%\quad\big(resp.(iT,R+iS,0,\cdots,0 )\big)
	.\label{v-rsting} \end{equation}
Thus by \eqref{gA} we have
\begin{eqnarray}
	f(v)&=&f(vO(v))=f\big((R+iS,iT,0,\cdots,0)\big)	%\quad\big(resp.g\big((iT,R+iS,0,\cdots,0)\big)\big),
	\label{I}\nonumber\\
	&=&\sqrt{v\bar{v}'}f\left(\Big(\frac{R+iS}{\sqrt{v\bar{v}'}},\frac{iT}{\sqrt{v\bar{v}'}},0,\cdots,0\Big)\right)\label{fA-rst},
\end{eqnarray}
by \eqref{v-r}-\eqref{v-st} we have
\begin{eqnarray}
	R&=&\frac{1}{2}\sqrt{vv'+\overline{vv'}+2v\bar{v}'}, \\
	S&=&\frac{1}{2i}\frac{vv'-\overline{vv'}}{\sqrt{vv'+\overline{vv'}+2v\bar{v}'}},  \\
	\frac{iT}{\sqrt{v\bar{v}'}}&=&
	\frac{i}{\sqrt{v\bar{v}'}}\sqrt{\frac{(v\bar{v}')^2-|vv'|^2}{vv'+\overline{vv'}+2v\bar{v}'}}
	=i\sqrt{\frac{1-t^2}{t+\bar{t}+2}}, \label{v-T}\\
	\frac{R+iS}{\sqrt{v\bar{v}'}}
	&=&\frac{1}{\sqrt{v\bar{v}'}}\frac{v\bar{v}'+vv'}{\sqrt{vv'+\overline{vv'}+2v\bar{v}'}}
	=\frac{1+t}{\sqrt{t+\bar{t}+2}}, \label{v-R+iS}
\end{eqnarray}
where we denote $t=\frac{vv'}{v\bar{v}'}$, and \eqref{fA-rst} can be writen as
\begin{eqnarray}
	f(v)=\sqrt{v\bar{v}'}f\left(\Big(\frac{1+t}{\sqrt{t+\bar{t}+2}},
	i\frac{\sqrt{1-t^2}}{\sqrt{t+\bar{t}+2}},0,\cdots,0\Big)\right)\label{fA-rst.},
\end{eqnarray}
the singular points appear when $ t=-1, $ which is the case when $v$ is pure imaginary, contradict to our assumption $ R\neq0 $. When $t(v)=1,$ $ \sqrt{1-t^2} $ is not smooth, but this is the case when $ S=T=0, $ which also contradict to our assumption.

By \eqref{gB} $f(v)=f(e^{i\theta}v)$ for $ \forall\theta\in \mathbb{R},$ and for $vv'=re^{i\theta},$ we have
$$ve^{-\frac{1}{2}i\theta}(\bar{v}e^{-\frac{1}{2}i\theta})'=v\bar{v}',
ve^{-\frac{1}{2}i\theta}(ve^{-\frac{1}{2}i\theta})'=r^2=|vv'|,$$
thus %we can assume  $t=\frac{|vv'|}{v\bar{v}'}$, and
\eqref{fA-rst.} can be writen as
\begin{eqnarray}
	f(v)=\frac{\sqrt{v\bar{v}'}}{2}f\left(\Big(\sqrt{1+t},
	i\sqrt{1-t},0,\cdots,0\Big)\right),\quad   t=\sqrt{s}=\frac{|vv'|}{v\bar{v}'}.
	\label{fA-rst1}
\end{eqnarray}

%Notice $ t $ is not smooth with respect to v on $T_{0}^{1,0} \mathfrak{R}_{IV}\setminus\{0\}.$ In order to make $f$ smooth with respect to v on $T_{0}^{1,0} \mathfrak{R}_{IV}\setminus\{0\},$ we have to consider $f$ as smooth function with respect to $s:=\frac{1}{4}\left[(\sqrt{1+t})^2-(\sqrt{1-t})^2\right]^2=t^2=\frac{|vv'|^2}{(v\bar{v}')^2}.$

For convenience, we denote $ \sqrt{\phi(s)}=\frac{1}{2}f\left(\Big(\sqrt{1+t},
i\sqrt{1-t},0,\cdots,0\Big)\right) ,$ thus
\begin{eqnarray}
	f(v)=\sqrt{r\phi(s)},\quad r=v\bar{v}', s(v)=\frac{|vv'|^2}{r^2}, 0\neq v\in \mathbb{C}^n \cong T_{0}^{1,0} \mathfrak{R}_{IV}.\label{g.}
\end{eqnarray}

Next we consider the condition of $f:\mathbb{C}^n\rightarrow [0,+\infty)$ to be  $\mbox{Iso}_0(\mathfrak{R}_{IV})$-invariant strongly pseudoconvex and strongly convex complex Minkowski norm, for any $f$ of the following form
	\begin{equation}
		f(\xi)=\sqrt{r\phi(s)},\label{cmn}
	\end{equation}
\begin{remark}
The invariant metrics considered in this subsection were constructed in \cite{GZ}. For the purpose of preserving the form of the original 2022 manuscript, no modifications have been made to the proof in this section.
\end{remark}

\section*{Acknowledgement}
This work was completed while the author was a PhD student at the School of Mathematical Sciences, Xiamen University.

\end{document}